\documentclass[10pt,twoside,reqno,a4paper]{amsproc}
\usepackage[utf8]{inputenc}
\usepackage{courier}
\usepackage[T1]{fontenc}

\usepackage[scaled]{helvet}
\usepackage[mathbf]{euler}
\usepackage{topformat}
\usepackage[driver=pdftex,margin=3cm,heightrounded=true,centering]{geometry}
\usepackage[backend=biber, style=alphabetic, url=false]{biblatex}
\RequirePackage{doi}

\usepackage{hyperref}
\usepackage{enumerate}
\usepackage{amssymb}
\usepackage{amsopn}
\usepackage{amsmath}
\usepackage{mathtools}
\usepackage{tikz}
\usepackage{tikz-cd}
\usepackage{ifthen}
\usepackage{graphics}
\usepackage{amsthm}
\usepackage{topthm}
\author{Thorben Kastenholz}
\thanks{Thorben Kastenholz was supported by the DFG (German Research
  Foundation) – SPP 2026, Geometry at Infinity - Project 73: Geometric Chern
  characters in p-adic equivariant K-theory}
\date{\today}
\title{Symplectic Groups, Mapping Class Groups and the Stability of Bounded
Cohomology}
\address{University of G\"ottingen, Bunsenstrasse 3-5, 37073 G\"ottingen,
Germany}
\email{thorben.kastenholz@mathematik.uni-goettingen.de}
\begin{document}
\newcommand{\introduce}[1]
  {\textbf{#1}}
\newcommand{\tk}[1]{\todo[size=\tiny,color=green!40]{TK: #1}}
\newcommand\blfootnote[1]{%
  \begingroup
  \renewcommand\thefootnote{}\footnote{#1}%
  \addtocounter{footnote}{-1}%
  \endgroup
}

\newcommand{\apply}[2]
  {{#1}\!\left({#2}\right)}
\newcommand{\at}[2]
  {\left.{#1}\right\rvert_{#2}}
\newcommand{\Identity}%
  {\mathrm{Id}}
\newcommand{\NaturalNumbers}%
  {\mathbf{N}}
\newcommand{\Integers}%
  {\mathbf{Z}}
\newcommand{\Rationals}%
  {\mathbf{Q}}
\newcommand{\Reals}%
  {\mathbf{R}}
  \newcommand{\ComplexNumbers}%
  {\mathbf{C}}
\newcommand{\AbstractProjection}[1] 
  {p_{#1}}
\newcommand{\RealPart}[1]
  {\apply{\operatorname{Re}}{#1}}
\newcommand{\ImaginaryPart}[1]
  {\apply{\operatorname{Im}}{#1}}
\newcommand{\Floor}[1]
  {\left \lfloor #1 \right \rfloor}
\newcommand{\Hom}[2]
  {\apply{\mathrm{Hom}}{#1,#2}}
\newcommand{\Norm}[1]
  {\left|\left|#1\right|\right|}
\newcommand{\AbsoluteValue}[1]
  {\left|#1\right|}
\newcommand{\colim}%
  {\mathrm{colim}}

\newcommand{\Manifold}%
  {M}
  \newcommand{\ManifoldAlternative}%
  {N}
 \newcommand{\ManifoldAuxiliary}%
  {K}
\newcommand{\NullBordism}%
  {W}
\newcommand{\Bordism}
  {P}
 \newcommand{\ManifoldFiber}%
 {M}
 \newcommand{\ManifoldTotal}%
 {E}
  \newcommand{\ManifoldBase}%
  {B}
\newcommand{\SmoothMap}%
  {\phi}
  \newcommand{\MorseFunction}%
  {f}
\newcommand{\Diffeomorphism}%
  {\Phi}
\newcommand{\Dimension}%
  {d}
  \newcommand{\HalfDimension}%
  {n}
\newcommand{\FundamentalClass}[1]
  {\left[#1\right]}
\newcommand{\Interval}%
  {I}
\newcommand{\Ball}[1]
  {D^{#1}}
\newcommand{\Sphere}[1] 
  {S^{#1}}
\newcommand{\Torus}[1]
  {T^{#1}}
\newcommand{\SimplicialVolume}[1]
  {\lvert \lvert #1 \rvert \rvert}
\newcommand{\ellone}%
  {\ell_{1}}
\newcommand{\Boundary}[1]
  {\partial #1}
\newcommand{\ComplexOfEmbeddings}[1]
  {\apply{K}{#1}}
\newcommand{\GenusOf}[1]
  {\apply{\Genus}{#1}}
\newcommand{\StableGenusOf}[1]
  {\apply{\overline{\Genus}}{#1}}
\newcommand{\Surface}[1]
  {\Sigma_{#1}}
\newcommand{\CurveComplex}[1]
  {\mathcal{C}_{#1}}
\newcommand{\MappingClassGroup}[1]
  {\apply{\mathcal{M}}{#1}}
\newcommand{\SignatureClass}[2]
  {\tau^{#1}_{#2}}
\newcommand{\TorelliGroup}[1]
  {\mathcal{T}_{#1}}

\newcommand{\Diff}[1]
  {\mathrm{Diff}\!\left(#1\right)}
\newcommand{\DiffBoundary}[1]
  {\apply{\mathrm{Diff}_\partial}{#1}}
\newcommand{\DiffGroup}[1]
  {\mathrm{Diff}^{B\Group}\!\left(#1\right)}
\newcommand{\DiffZero}[1]
  {\mathrm{Diff}_0\!\left(#1\right)}
\newcommand{\DiffOne}[1]
  {\widetilde{\mathrm{Diff}}_0\!\left(#1\right)}
\newcommand{\HomeoGroup}[1]
  {\apply{\mathrm{Homeo}}{#1}}
\newcommand{\HomeoCompactGroup}[1]
  {\apply{\mathrm{Homeo}_{c}}{#1}}
\newcommand{\HomeoLowerGroup}[1]
  {\apply{\mathrm{Homeo}^{\geq}}{#1}}

\newcommand{\FiberingSpace}%
  {E}
\newcommand{\FiberingProjektion}[1] 
  {\pi_{#1}}
\newcommand{\Fiber}%
  {F}
\newcommand{\FiberDimension}%
  {f}
\newcommand{\Base}%
  {B}
\newcommand{\ClutchingFunction}[1] 
  {\varphi_{#1}}

\newcommand{\Group}%
  {\Gamma}
\newcommand{\GroupElement}%
  {g}
\newcommand{\Genus}%
  {g}
\newcommand{\QuadraticModule}%
  {\mathbf{M}}
\newcommand{\WittIndex}[1]
  {\apply{\Genus}{#1}}
\newcommand{\StableWittIndex}[1]
  {\apply{\overline{\Genus}}{#1}}
\newcommand{\ComplexOfHyperbolicInclusions}[1]
  {\apply{K^{a}}{#1}}
\newcommand{\ChainContraction}[1]
  {H_{#1}}
\newcommand{\SymplecticGroup}[2]
  {\apply{\mathrm{Sp}_{#1}}{#2}}

\newcommand{\HomologyClass}%
  {\alpha}
\newcommand{\HomologyOfSpaceObject}[3]
  {\apply{H_{#1}}{#2 ; #3}}
\newcommand{\CohomologyOfSpaceObject}[3]
  {\apply{H^{#1}}{#2 ; #3}}
\newcommand{\BoundedCohomologyOfSpaceObject}[3]
  {\apply{H^{#1}_{\text{b}}}{#2 ; #3}}
\newcommand{\BoundedCohomologyOfSimplicialObject}[3]
  {\apply{H^{#1}_{\text{b, s}}}{#2 ; #3}}
\newcommand{\HomologyOfSpaceMorphism}[1]
  {{#1}_{\ast}}
\newcommand{\CohomologyOfSpaceMorphism}[1]
  {{#1}^{\ast}}
\newcommand{\HomologyOfGroupObject}[3]
  {\apply{H_{#1}}{#2; #3}}
\newcommand{\HomologyOfGroupMorphism}[1]
  {{#1}_{\ast}}
\newcommand{\HomologyOfSpacePairObject}[3]
  {\apply{H_{#1}}{{#2},{#3}}}
\newcommand{\Multiple}%
  {\lambda}

\newcommand{\TopologicalSpace}%
  {X}
\newcommand{\Point}%
  {\ast}
\newcommand{\Loop}%
  {\gamma}
\newcommand{\ContinuousMap}%
  {f}
  \newcommand{\ContinuousMapALT}%
  {g}
 \newcommand{\maps}%
  {\ensuremath{\text{maps}}}
\newcommand{\HomotopyGroupOfObject}[3]
  {\apply{\pi_{#1}}{{#2},{#3}}}
\newcommand{\HomotopyGroupOfPairObject}[4]
  {\apply{\pi_{#1}}{{#2},{#3},{#4}}}
\newcommand{\HomotopyGroupMorphism}[1] 
  {{#1}_{\ast}}
\newcommand{\EMSpace}[2]
  {\apply{K}{{#1},{#2}}}
\newcommand{\ClassifyingSpace}[1] 
  {B#1}
\newcommand{\UniversalCovering}[1] 
  {\widetilde{#1}}
\newcommand{\UniversalCoveringMap}[1] 
  {\widetilde{#1}}
\newcommand{\Signature}[1]
  {\apply{\sigma}{#1}}

\newcommand{\SimplicialComplex}
  {X}
\newcommand{\AuxSimplicialComplex}
  {K}
\newcommand{\Subcomplex}
  {Y}
\newcommand{\AuxSubcomplex}
  {L}
\newcommand{\Simplex}[1]
  {\sigma_{#1}}
\newcommand{\Link}[2]
  {\apply{\text{Lk}_{#1}}{#2}}
\newcommand{\Star}[2]
  {\apply{\text{St}_{#1}}{#2}}
\newcommand{\BoundaryIndexSimplex}[2]
  {\apply{\partial_{#1}}{#2}}
\newcommand{\BoundarySimplex}%
  {\partial}
\newcommand{\StandardSimplex}[1]
  {\Delta_{#1}}
\newcommand{\vertex}
  {v}
\newcommand{\GeometricRealization}[1]
  {\left\lvert #1 \right\rvert}
\newcommand{\BoundHomotopy}
  {N}
\newcommand{\Horn}[2]
  {\Lambda^{#1}_{#2}}

\begin{abstract}
  Mapping class groups satisfy cohomological stability. In this note we show how
  results by Bestvina and Fujiwara imply that their bounded cohomology does not
  stabilize, additionally we show that stabily polynomials in the
  Mumford-Morita-Miller classes are unbounded i.e. their norm tends to
  infinity as one increases the genus.

  While the bounded cohomology of the symplectic group does stabilize, we show
  that it does not stabilize via isometries in degree $2$.
  In order to establish this we calculate the norm of the signature class in
  $\SymplecticGroup{2h}{\Reals}$ and estimate the norm of the integral
  signature class.
\end{abstract}

\maketitle
\section{Introduction}
Homological stability refers to a phenomenon, where the group homology of an increasing
sequence of groups $G_0 \subset G_1 \subset \ldots$ eventually
stabilizes i.e. the inclusions induce isomorphisms in more and more degrees.
The two most important examples for our purposes are given by symplectic groups
$\SymplecticGroup{R}{2h}$ and mapping class groups of surfaces. More precisely,
consider the mapping class group $\MappingClassGroup{\Surface{g,b}}$ of a
surface of genus $g$ with $b$ boundary components, then there are
maps
\begin{align*}
  \MappingClassGroup{\Surface{g,b}}
  &
  \to
  \MappingClassGroup{\Surface{g + 1,b}}
  \\
  \MappingClassGroup{\Surface{g,b}}
  &
  \to
  \MappingClassGroup{\Surface{g,b+1}}
  \\
  \MappingClassGroup{\Surface{g,b}}
  &
  \to
  \MappingClassGroup{\Surface{g,b-1}}
\end{align*}
coming from the various inclusions of surfaces and they all induce
isomorphisms in homology in a range of degrees that increases with $g$.

Classically homological stability is proven by constructing a sequence of
simplicial complex $X_i$ on which $G_i$ acts in a particular fashion and which
have to be acyclic in an increasing range of degrees. This approach goes back
to Quillen.
Recently this approach has been adapted by Hartnick and De La Cruz Mengual in
\cite{HartnickDeLaCruzQuillen} to yield a similar setup that proves stability
for
bounded cohomology. The most important change lies in the fact that the
simplicial complexes have to have vanishing simplicial bounded cohomology in an
increasing range of degrees. Many of the necessary tools to translate a
classical homological stability proof to the realm of bounded cohomology have
been established in \cite{KastenholzSroka}.
There, one of the key observations is that a simplicial complex with vanishing
first bounded cohomology has to have finite diameter (See also Example~something
in \cite{KastenholzSroka}) and that an acyclic complex has to satisfy a higher
dimensional analogue of having finite diameter in order to have vanishing
simplicial bounded cohomology in range.

The complexes in question for the mapping class groups are the so called curve
complexes. While these are highly connected, they have infinite diameter.
Furthermore they are hyperbolic as shown in \cite{MazurMinski}. These two
facts allow one to construct quasimorphisms that will lie in the kernel of the
corresponding inclusions in bounded cohomology or in other words show that the
map induced by the inclusion
\[
  \BoundedCohomologyOfSpaceObject{2}{\MappingClassGroup{\Surface{g,b}}}{\Reals}
  \to
  \BoundedCohomologyOfSpaceObject{2}{\MappingClassGroup{\Surface{g+1,b}}}{\Reals}
\]
is never injective. This will be carried out in Section~\ref{scn:MCG}.
This instability is closely related to one of Kirby's Problems
(Problem~2.18 in \cite{KirbyList}), for a given $h$, what is the
minimal $g$ such that there exists a fiber bundle
\[
  \Surface{h}\to E \to \Surface{g}
\]
with $E$ having non-zero or fixed signature? Indeed by stabilizing this
question as follows: Let $\apply{g_h}{n}$ denote the minimal genus $g$ such
that there
exists a surface bundle $\Surface{h}\to E\to \Surface{g}$ with signature $4n$
(the signature is always a multiple of $4$ by \cite{MayerSignature}). Now let
$G_h$ denote the limit of $\frac{\apply{g_h}{n}}{n}$, then one can show that
$G_h$ corresponds to the $\ellone$-norm of a generator of
$\HomologyOfGroupObject{2}{\MappingClassGroup{\Surface{h}}}{\Integers}$. By
duality computing this norm is equivalent to computing the norm of the
generator of the second group cohomology of the mapping class group.
Hence the aforementioned instability encapsulates how $G_h$ changes, when
changing $h$. Furthermore by now there exist enough examples to show that $G_h$
tends to zero as $h$ increases (see for example \cite{Monden}), which also
shows that the
norm of the corresponding generator of the second group cohomology tends to
infinity. Finally the stable cohomology of mapping class groups is a polynomial
ring generated by the so called Mumford-Morita-Miller-Classes. These are
characteristic
classes that are defined by fiber integrating powers of the Euler class
of the vertical tangent bundle of a surface bundle. In
Proposition~\ref{prp:MMMClassesUnstable}, we show that the norm of these
classes also tends to infinity as the genus increases. In particular the
MMM-classes on the stable mapping class group cannot be bounded (See
Corollary~\ref{cor:StableMCG}). This is
particularly interesting, since it was shown in \cite{MoritaBoundedness}, that
the odd MMM-classes are bounded and the recent preprint \cite{Hamenstaedt},
claims that all MMM-classes are bounded.

Via the action of the mapping class group on the second homology of the
surface one obtains a close relationship between mapping class groups and
symplectic groups. In particular the aforementioned map induces an isomorphism
on second group cohomology. It was shown in \cite{HartnickDeLaCruz} that the
bounded cohomology of symplectic groups stabilizes. This difference between
mapping class groups and symplectic groups together with their relationship
yields some interesting properties of the stability of the bounded cohomology
of symplectic groups.
First and foremost using this one can show that the stabilization maps
between integral as well as real symplectic groups while inducing isomorphism
fail to induce isometries on second bounded cohomology. In \cite{DCMThesis} it
was already shown via explicit computations that the third bounded cohomology
of symplectic groups does not stabilize isometrically (their result is about
continuous
bounded cohomology, but using Corollary~1.4 in \cite{MonodLattice}, this can be
translated to discrete bounded cohomology of symplectic groups of the Gaussian
integers).
The herein presented example is particularly interesting since this occurs in
degree $2$, where bounded cohomology is in
fact a Banach space (see \cite{MatsumotoMorita}). More precisely one can show
that the maps $\SymplecticGroup{2h}{R}\to \SymplecticGroup{2(h+k)}{R}$ have
arbitrarily small norm as $k$ tends to infinity for $R$ either the reals or
integers.

Furthermore as a consequence one obtains that the bounded cohomology of the
stable symplectic group $\SymplecticGroup{\infty}{R}$  (i.e.
the colimit of the aforementioned sequence) does not surject onto the stable
bounded cohomology of the finite dimensional symplectic groups (indeed this
would imply
that the aforementioned inclusions do not have arbitrarily small norm). Here
one should note the relationship to the vanishing modulus introduced in
Section~4 in \cite{BinateGroups}, where the relationship between the vanishing
of bounded cohomology and colimits was investigated. All in all this shows that
the
stability of bounded cohomology seems to be more delicate than ordinary
(co)homological stability.

In particular one can ask the general question:
\begin{question}
  If the bounded cohomology of a
  sequence of groups $\Group_i$ does stabilize, is there an injection
  \[
    \BoundedCohomologyOfSimplicialObject{*}{\colim \Group_i}{\Reals}
    \to
    \colim
    \BoundedCohomologyOfSpaceObject{*}{\Group_i}{\Reals}
  \]
  that hits all classes with overall bounded norm?
\end{question}
Since a lot of the usage of classical homological stability results comes from
the stable group $\colim \Group_i$ behaving more nicely than $\Group_i$,
answering this question would
probably allow one to get a better understanding of the bounded cohomology of
groups that satisfy stability of bounded cohomology.

One the other side, the proof of the instability of the bounded cohomology of
mapping class groups relied on the occuring curve complexes having infinite
diameter so the obvious question becomes:
\begin{question}
    If the complexes occuring in a
  classical proof of homological stability have infinite diameter, can one
  conclude that the second bounded cohomology of these groups does not
  stabilize?
\end{question}
Stated that way, the question is probably too naive. In particular the proof
not only relies on the curve complexes being hyperbolic, but also on the fact
that the curve complex of $\Surface{g,1}$ injects into the curve complex
$\Surface{g+1,1}$ such that the image has finite diameter.

Lastly the author does not claim any originalty. This note represents the
application of known results in the context of stability of bounded cohomology.
Finally the author would like to thank Robin Sroka and Mladen Bestvina for
helpful discussions. Lastly the author thanks Mark Pedron for providing the
argument in the proof of Proposition~\ref{prp:MMMClassesUnstable}.
\section{Mapping Class Groups}
\label{scn:MCG}
The well-studied curve complex $\CurveComplex{h}$ of a surface $\Surface{h}$
is defined as follows:
The vertices are given by isotopy classes of simple closed curves $\gamma$ on
$\Surface{h}$ and a collection of curves $(\gamma_0,\ldots,\gamma_p)$ forms a
$p$-simplex if the isotopy classes can be realized with pairwise disjoint
images. This complex and its many variatons have an intricate relationship with
the mapping class group. In particular they were used by Harer in
\cite{HarerHomologicalStability} to establish homological stability for
mapping class groups. Bestvina and Fujiwara used the curve complex to construct
many non-trivial quasimorphisms on $\MappingClassGroup{\Surface{h}}$ and we
will use these quasimorphisms to show that the inclusion
\[
  \MappingClassGroup{\Surface{h,1}} \to \MappingClassGroup{\Surface{h+1,1}}
\]
will never be injective. The mapping class group acts simplicially on the
curve complex via postcomposition of the curves.
Let $\apply{d}{x,y}$ denote the simplicial distance function on the one
skeleton of the curve complex denoted by $\CurveComplex{h}^{(1)}$. Given two
paths $w$ and $\alpha$ of finite length in $\CurveComplex{h}^{(1)}$, we define
$
  \AbsoluteValue{\alpha}_{w}
$
to be the maximal number of non-overlapping copies of $w$ in $\alpha$. Now
for every $0 < W < \AbsoluteValue{w}$, Bestvina and Fujiwara define
\[
  \apply{c_{w,W}}{x,y}
  =
  \apply{d}{x,y}
  -
  \inf_\alpha
  \left(
    \AbsoluteValue{\alpha}
    -
    W\AbsoluteValue{\alpha}_{w}
  \right)
\]
where $\alpha$ ranges over all paths from $x$ to $y$. Now they fix a $W$ that
is big compared to some constant depending only on the curve complex and
define
\begin{align*}
  h_{w,x_0}
  \colon
  \MappingClassGroup{\Surface{h}}
  &
  \to
  \Reals
  \\
  g
  &
  \mapsto
  \apply{c_{w}}{x_0,\apply{g}{x_0}}
  -
  \apply{c_{w^{-1}}}{x_0,\apply{g}{x_0}}
\end{align*}
where $w^{-1}$ denotes $w$ with the opposite orientation. They then proceed
to show that for any given $x_0$, these define quasimorphisms on the mapping
class group and that under the correspondence between quasimorphisms and
second bounded cohomology these $h_w$ span an infinite dimensional subspace
of
$E\BoundedCohomologyOfSpaceObject{2}{\MappingClassGroup{h}}{\Reals}$.

Now consider the sequence of inclusions
\[
  \Surface{g,1} \subset \Surface{g+1,1} \subset \Surface{g+1}
\]
which induce inclusions
\[
  \MappingClassGroup{\Surface{g,1}}
  \subset
  \MappingClassGroup{\Surface{g+1,1}}
  \subset
  \MappingClassGroup{\Surface{g+1}}
\]
that in turn induce the isomorphisms in Harer's stability result.
Let $x_0$ denote a curve in $\Surface{g+1,1}$ that is disjoint from
$\Surface{g,1}$, then $\MappingClassGroup{\Surface{g,1}}$ fixes this element.
Hence for any $w$ we have that $h_{w,x_0}$  is actually zero when restricted
to
$
  \MappingClassGroup{\Surface{g,1}}
$%
. In particular the bounded cohomology of $\MappingClassGroup{\Surface{g,1}}$
does not stablizie with respect to $g$ since the map induced by the inclusion
in degree two will always have a kernel.

\subsection{MMM-Classes are stably unbounded}
As mentioned in the introduction, mapping class groups satisfy homological
stability. Moreover in \cite{MadsenWeiss}, the stable cohomology of mapping
class groups was computed and it is a polymial ring generated by the so called
Mumford-Morita-Miller-Classes. These are defined as follows:
Let $\FiberingProjektion{\FiberingSpace} \colon \FiberingSpace\to \Base$ denote
an oriented surface bundle and $T_v
\FiberingSpace$ its vertical tangent bundle. Let $\apply{e}{T_v
\FiberingSpace}$ denote the Euler class of the vertical tangent bundle. Now we
define the $k$-th MMM-class via
\[
  \apply{\kappa_k}{\FiberingSpace}
  \coloneqq
  \apply
    {\left(\FiberingProjektion{\FiberingSpace}\right)_{!}}
    {\apply{e^{k+1}}{T_v \FiberingSpace}}\text{.}
\]
\begin{lemma}
  Let $f \colon \hat{\FiberingSpace}\to \FiberingSpace$ denote a degree $d$
  covering map
  that is fiberwise a covering of the bundle
  $
    \FiberingProjektion{\FiberingSpace}
    \colon
    \FiberingSpace
    \to
    \Base
  $,
  then
  $
    \apply{\kappa_k}{\hat{\FiberingSpace}}
    =
    d \apply{\kappa_k}{\FiberingSpace}
  $%
  .
\end{lemma}
\begin{proof}
  Let $\FiberingProjektion{\hat{\FiberingSpace}}$ denote the composition of $f$
  and $\FiberingProjektion{\FiberingSpace}$.
  The fiber integration of a composition is the composition of the fiber
  integrations. Since the fiber integration of a covering map is given by the
  transfer $tf$, we have
  \begin{align*}
    \apply{\kappa_k}{\hat{\FiberingSpace}}
    &
    =
    \apply
      {(\FiberingProjektion{\hat{\FiberingSpace}})_{!}}
      {\apply{e^{k+1}}{T_v \hat{\FiberingSpace}}}\\
    &
    =
    \apply
      {(\FiberingProjektion{\FiberingSpace}\circ f)_{!}}
      {\apply{e^{k+1}}{f^*T_v \FiberingSpace}}\\
    &
    =
    \apply
      {(\FiberingProjektion{\FiberingSpace})_{!}\circ tf}
      {\apply{f^*}{\apply{e^{k+1}}{T_v \FiberingSpace}}}\\
    &
    =
    \apply
      {(\FiberingProjektion{\FiberingSpace})_!}
      {d \apply{e^{k+1}}{T_v \FiberingSpace}}\\
    &
    =
    d \apply{\kappa_k}{\FiberingSpace}
  \end{align*}
\end{proof}
With this lemma at hand we can prove that stabily all polynomials in the
Mumford-Morita-Miller classes are unbounded.
\begin{proposition}
\label{prp:MMMClassesUnstable}
  Let
  $
    P
    \in
    \lim_g
    \CohomologyOfSpaceObject{k}{\MappingClassGroup{\Surface{g}}}{\Reals}
  $
  denote a stable cohomology class, then $\Norm{P_g}$ tends to infinity as $g$
  increases, where $P_g$ denotes the restriction of $P$ to
  $\MappingClassGroup{\Surface{g}}$.
\end{proposition}
\begin{proof}
  By \cite{MadsenWeiss}, we have
  $
    \lim
    \CohomologyOfSpaceObject{*}{\MappingClassGroup{\Surface{g}}}{\Reals}
    \cong
    \Reals[(\kappa_i)_{i\in \NaturalNumbers_{>0}}]
  $
  in particular $P$ is a polynomial in $\kappa_k$, which is homogenous with
  respect to the cohomology grading. Let $\hat{P} = \kappa_i^d$ denote one of
  its non-zero summands, then there exists a $k$-dimensional manifold $\Base$,
  together with
  a surface bundle
  $
    \FiberingProjektion{\FiberingSpace}
    \colon
    \FiberingSpace
    \to
    \Base
  $
  with fiber genus $g$ for $g$ large enough, such that
  $\apply{\hat{P}_g}{\FiberingSpace}$ is $\lambda$, but all other summands of
  $P_g$ evaluate to zero on $\FiberingSpace$.

  The kernel of the composition
  $
    \HomotopyGroupOfObject{1}{\Surface{g}}{\ast}
    \to
    \HomologyOfSpaceObject{1}{\Surface{g}}{\Integers}
    \to
    \HomologyOfSpaceObject{1}{\Surface{g}}{\Integers/n\Integers}
  $
  is invariant under the action of the diffeomorphism group of $\Surface{g}$ on
  the fundamental group, hence it defines a fiberwise finite covering
  $
    \hat{\FiberingSpace}
    \to
    \FiberingSpace
    \to
    \Base
  $%
  . Let us define the composition of these two maps again by
  $\FiberingProjektion{\hat{\FiberingSpace}}$; this map is again a surface
  bundle with fiber genus $g_n = n^{2g}(g-1) - 1$.
  Now we have
  \[
    (n^{2g})^d \lambda
    =
    \langle
      \apply{\hat{P}_{g_n}}{\hat{\FiberingSpace}}
      ,
      \FundamentalClass{\Base}
    \rangle
    =
    \langle
      \apply{P_{g_n}}{\hat{\FiberingSpace}}
      ,
      \FundamentalClass{\Base}
    \rangle
    \leq
    \Norm{P_{g_n}}
    \Norm{B}
  \]
  and hence $\Norm{P_{g_n}}$ tends to infinity as $n$ tends to infinity.
\end{proof}
This immediately yields the following corollary:
\begin{corollary}
\label{cor:StableMCG}
  Let $\MappingClassGroup{\infty}$ denote the stable mapping class group i.e.
  $
    \colim_g \MappingClassGroup{\Surface{g,1}}
  $%
  , then the comparison map
  \[
    \BoundedCohomologyOfSpaceObject{*}{\MappingClassGroup{\infty}}{\Reals}
    \to
    \CohomologyOfSpaceObject{*}{\MappingClassGroup{\infty}}{\Reals}
  \]
  is the zero map.
\end{corollary}
\section{Symplectic Groups}
The second group cohomology with integer coefficients of the symplectic groups
$\SymplecticGroup{2h}{\Integers}$ and $\SymplecticGroup{2h}{\Reals}$ is
generated by the signature class $\SignatureClass{R}{2h}$, where $R$ is either
$\Integers$ or $\Reals$, defined as follows:
Since the symplectic group is perfect one has that its first homology
vanishes, hence
$\CohomologyOfSpaceObject{2}{\SymplecticGroup{2h}{R}}{\Integers}$ is
isomorphic to
$\Hom{\HomologyOfGroupObject{2}{\SymplecticGroup{2h}{R}}{\Integers}}{\Integers}$.
Now one can represent every second homology class by a map
$
  f
  \colon
  \Surface{g}
  \to
  \ClassifyingSpace{\SymplecticGroup{2h}{R}}
$
which corresponds to a twisted coefficient system $E$ on $\Surface{g}$. One can
check that this induces a symmetric bilinear pairing
\[
  \CohomologyOfSpaceObject{1}{\Surface{g}}{E}
  \times
  \CohomologyOfSpaceObject{1}{\Surface{g}}{E}
  \to
  \CohomologyOfSpaceObject{2}{\Surface{g}}{E}
\]
and
$
  \apply
    {\SignatureClass{R}{2h}}
    {
      \apply
        {f}
        {\FundamentalClass{\Surface{g}}}
    }
$
is defined to be the signature of this symmetric bilinear form.
We are interested in the norm of these classes.
Since there is an inclusion
$
  \SymplecticGroup{2h}{\Integers}
  \to
  \SymplecticGroup{2h}{\Reals}
$
we have
$
  \Norm{\SignatureClass{\Integers}{2h}}
  \leq
  \Norm{\SignatureClass{\Reals}{2h}}
$%
.
Meyer in \cite{MayerSignature} and later Turaev in \cite{Turaev} gave explicit
cocycles
for these classes.
These cocycles yield
\[
  \Norm{\SignatureClass{R}{2h}} \leq 2h
\]%
but one can do better in the real case: In \cite{LionVergne} they construct a
cocycle
representative of $\SignatureClass{\Reals}{2h}$ in terms of Lagrangian
subspaces of $\Reals^{2h}$ and the Maslov index that has norm bounded
by $h$ (See also Section~5 in \cite{Turaev}).
Now note that Milnor has computed $\Norm{\SignatureClass{\Reals}{2}}$ to be
$1$ by giving examples of flat $\SymplecticGroup{\Reals}{2}$-bundles over a
surface of genus $g$ with Euler class $g-1$.
Since the signature class is equal to four times the Euler class (again see
\cite{Turaev}), these are examples where the signature is maximal. Since the
signature is additive with respect to direct sums of flat
$\SymplecticGroup{\Reals}{2}$-bundles, we obtain that
$\Norm{\SignatureClass{\Reals}{2h}}$ is indeed $h$.

Note that by \cite{Charney} the symplectic groups satisfy
homological stability. Additionally by \cite{HartnickDeLaCruz} and
\cite{KastenholzSroka} their bounded
cohomology stabilizes as well. Hence we have diagrams:
\[
  \begin{tikzcd}
    \BoundedCohomologyOfSpaceObject{2}{\SymplecticGroup{\Reals}{2h}}{\Reals}
      \ar[r]
      \ar[d]
    &
    \BoundedCohomologyOfSpaceObject{2}{\SymplecticGroup{\Reals}{2h+2k}}{\Reals}
      \ar[d]
    \\
    \CohomologyOfSpaceObject{2}{\SymplecticGroup{\Reals}{2h}}{\Reals}
      \ar[r]
    &
    \CohomologyOfSpaceObject{2}{\SymplecticGroup{\Reals}{2h+2k}}{\Reals}
  \end{tikzcd}
\]
where both horizontal arrows are isomorphisms. In particular if the upper
isomorphism would be an isometry, then the lower map would be an isometry as
well since the norm on ordinary cohomology is the quotient norm with respect to
the vertical maps.
Hence while the inclusion
\[
  \SymplecticGroup{\Reals}{2h}
  \to
  \SymplecticGroup{\Reals}{2h + 2k}
\]
induces isomorphisms on second bounded cohomology for $h$ big enough, these
fail to be isometries.

\subsection{Relation to the mapping class group}
Note that the action of the mapping class group on the first cohomology of the
corresponding surface yields a surjection
\[
  \Phi
  \colon
  \MappingClassGroup{\Surface{h}}
  \to
  \SymplecticGroup{\Integers}{2h}
\]
and given a surface bundle over a surface $\Surface{h} \to E \to \Surface{g}$,
the signature of $E$ is exactly given by
$
  \apply
    {\apply
      {\CohomologyOfSpaceMorphism{(\psi_E \circ \Phi)}}
      {\SignatureClass{\Integers}{2h}}
    }
    {\FundamentalClass{\Surface{g}}}
$%
, where $\psi_E$ denotes a classifying map of $E$.
This implies that
\[
  \Norm{
    \apply
      {\CohomologyOfSpaceMorphism{\Phi}}
      {\SignatureClass{\Integers}{2h}}
  }
  \leq
  \Norm{\SignatureClass{\Integers}{2h}}
\]

Albeit we will not use this, there is more to be said about this morphism:
Additionally $\Phi$ induces an isomorphism in second group cohomology and by
by Theorem~2.14 in \cite{Bucher} $\Phi$ induces
an injective isometry in bounded cohomology, but by
the observations in Section~\ref{scn:MCG} fails to be an
isomorphism. Additionally results by Hamenst\"adt in \cite{Hamenstaedt} and
before that by Kotschik \cite{Kotschick} show that the norm of
$
  \apply
    {\CohomologyOfSpaceMorphism{\Phi}}
    {\SignatureClass{\Integers}{2h}}
$
is bounded by $\frac{1}{2}(h-1)$ and $\frac{1}{3}(h-1)$ respectively.

In the paper \cite{EndoKorkmazKotschickOzbagciStipsicz}, for any
$h\geq3$ they construct a surface bundle $E$ with fiber genus $h$ over a
surface of
genus $9$ with signature $4\frac{h-2}{3}$. Together with the previous
inequality we obtain
\[
  4\frac{h-2}{3}
  =
    \apply
      {\apply
        {\CohomologyOfSpaceMorphism{(\psi_E \circ \Phi)}}
        {\SignatureClass{\Integers}{2h}}
      }
      {\FundamentalClass{\Surface{9}}}
  \leq
  \Norm{
      \apply
      {\CohomologyOfSpaceMorphism{\Phi}}
      {\SignatureClass{\Integers}{2h}}
  }
  \Norm{\Surface{9}}
  \leq
  32
  \Norm{\SignatureClass{\Integers}{2h}}
\]
showing that the norm of the integral signature class tends to infinity with
increasing $h$ as well.

Hence the norm of $\SignatureClass{\Integers}{2h}$ lies somewhere inbetween
$\frac{h-2}{24}$ and $h$. Showing that the inclusions between integral
symplectic groups while inducing isomorphisms on second bounded cohomology do
not induce isometries as well.

\emergencystretch=10em
\printbibliography
\end{document}